\documentclass[12pt,thmsa]{article}
\usepackage{amsfonts}
\usepackage{sw20bams}

\input tcilatex
\begin{document}

\author{Steven Finch}
\title{Powers of Euler's $q$-Series}
\date{January 19, 2007}
\maketitle

\begin{abstract}
What are the asymptotic moments of coefficients obtained when expanding $%
\tprod\nolimits_{m=1}^\infty (1-q^m)^\ell $ in series? A few examples are
given, as well as a new multiplicative representation for coefficients when $%
\ell =10$ and $\ell =14$.
\end{abstract}

\footnotetext{
Copyright \copyright\ 2007 by Steven R. Finch. All rights reserved.}Starting
with Euler's result: 
\[
\dprod\limits_{m=1}^\infty (1-q^m)=\dsum\limits_{k=-\infty }^\infty
(-1)^kq^{k(3k+1)/2} 
\]
and Jacobi's cubic analog: 
\[
\dprod\limits_{m=1}^\infty (1-q^m)^3=\dsum\limits_{k=0}^\infty
(-1)^k(2k+1)q^{k(k+1)/2}, 
\]
it is natural to wonder about the corresponding square and $\ell ^{\text{th}%
} $ powers for $\ell \geq 4$. We will consider the coefficients $%
\{a_n\}_{n=0}^\infty $ of the $q$-series expansion (on the right-hand side)
for fixed $\ell $. Estimates of the magnitude of a finite subsequence $%
\{a_n\}_{n\leq N}$ of coefficients include 
\[
\begin{array}{ccccc}
\dsum\limits_{n\leq N}|a_n|, &  & \dsum\limits_{n\leq N}a_n^2, &  & 
\max\limits_{n\leq N}\,|a_n|
\end{array}
\]
and we will analyze the asymptotics (when possible)\ of these as $%
N\rightarrow \infty $. Much of the material presented here is perhaps not
new.\ Our contribution (as in \cite{FS1, FS2}) is only to collect results in
one place.

\section{First Power}

Let the coefficients $\{a_n\}_{n=0}^\infty $ satisfy \cite{Euler, Andrws,
Sec1} 
\[
\dsum\limits_{n=0}^\infty a_nq^n=\dprod\limits_{m=1}^\infty (1-q^m). 
\]
Euler's pentagonal-number theorem gives that 
\[
a_n=\left\{ 
\begin{array}{lll}
(-1)^k &  & \text{if }24n+1=(6k+1)^2\text{ for some }-\infty <k<\infty , \\ 
0 &  & \text{otherwise}
\end{array}
\right. 
\]
and hence, by elementary considerations, 
\[
\dsum\limits_{n\leq N}|a_n|\sim \frac 13(24N)^{1/2} 
\]
as $N\rightarrow \infty $. Also $a_n=f(24n+1)$, where $f$ is multiplicative
and 
\[
f(p^r)=\left\{ 
\begin{array}{lll}
1 &  & \text{if }p\equiv 1,11\limfunc{mod}12\text{ and }r\text{ is even,} \\ 
(-1)^{r/2} &  & \text{if }p\equiv 5,7\limfunc{mod}12\text{ and }r\text{ is
even,} \\ 
0 &  & \text{otherwise}
\end{array}
\right. 
\]
for $p$ prime and $r\geq 0$. Of course, $\sum_{n\leq N}a_n$ oscillates
between $-1,0,1$ and never converges.

\section{Second Power}

Let the coefficients $\{a_n\}_{n=0}^\infty $ satisfy \cite{Rama, Gls1, Sec2} 
\[
\dsum\limits_{n=0}^\infty a_nq^n=\dprod\limits_{m=1}^\infty (1-q^m)^2. 
\]
No simple formula for $a_n$ is known, but $a_n=f(12n+1)$, where $f$ is
multiplicative and 
\[
f(p^r)=\left\{ 
\begin{array}{lll}
1 &  & \text{if }p\equiv 7,11\limfunc{mod}12\text{ and }r\text{ is even,} \\ 
(-1)^{r/2} &  & \text{if }p\equiv 5\limfunc{mod}12\text{ and }r\text{ is
even,} \\ 
r+1 &  & \text{if }p\equiv 1\limfunc{mod}12\text{ and }(-3)^{(p-1)/4}\equiv 1%
\limfunc{mod}p\text{,} \\ 
(-1)^r(r+1) &  & \text{if }p\equiv 1\limfunc{mod}12\text{ and }%
(-3)^{(p-1)/4}\equiv -1\limfunc{mod}p\text{,} \\ 
0 &  & \text{otherwise.}
\end{array}
\right. 
\]
We examine 
\begin{eqnarray*}
\dsum\limits_{n=1}^\infty \frac{|f(n)|}{n^s} &=&\dprod\limits\Sb p\equiv
5,7,11  \\ \limfunc{mod}12  \endSb \left( 1+\dsum\limits_{r=1}^\infty \frac
1{p^{2rs}}\right) \cdot \dprod\limits\Sb p\equiv 1  \\ \limfunc{mod}12 
\endSb \left( 1+\dsum\limits_{r=1}^\infty \frac{r+1}{p^{rs}}\right) \\
\ &=&\dprod\limits\Sb p\equiv 5,7,11  \\ \limfunc{mod}12  \endSb \left(
1-\frac 1{p^s+1}\right) \left( 1-\frac 1{p^s}\right) ^{-1}\cdot
\dprod\limits \Sb p\equiv 1  \\ \limfunc{mod}12  \endSb \left( 1+\frac
1{p^s-1}\right) \left( 1-\frac 1{p^s}\right) ^{-1} \\
\ &=&\left( 1-\frac 1{2^s}\right) ^{1/2}\left( 1-\frac 1{3^s}\right)
^{1/2}\dprod\limits\Sb p\equiv 5,7,11  \\ \limfunc{mod}12  \endSb \left(
1+\frac 1{p^s}\right) ^{-1}\left( 1-\frac 1{p^s}\right) ^{-1/2} \\
&&\ \ \cdot \dprod\limits\Sb p\equiv 1  \\ \limfunc{mod}12  \endSb \left(
1-\frac 1{p^s}\right) ^{-3/2}\cdot \zeta (s)^{1/2}
\end{eqnarray*}
and deduce, via the Selberg-Delange method \cite{Slbrg, Dlnge, Tnbm1}, that 
\[
\dsum\limits_{n\leq N}|a_n|\sim C\frac N{(\ln N)^{1/2}} 
\]
as $N\rightarrow \infty $, where \cite{Sebah} 
\begin{eqnarray*}
C &=&\frac{12}{\sqrt{3\pi }}\dprod\limits\Sb p\equiv 5,7,11  \\ \limfunc{mod}%
12  \endSb \left( 1+\frac 1p\right) ^{-1}\left( 1-\frac 1p\right)
^{-1/2}\cdot \dprod\limits\Sb p\equiv 1  \\ \limfunc{mod}12  \endSb \left(
1-\frac 1p\right) ^{-3/2} \\
\ &=&3.3215840614847482046694103....
\end{eqnarray*}
Also, 
\begin{eqnarray*}
\dsum\limits_{n=1}^\infty \frac{f(n)^2}{n^s} &=&\dprod\limits\Sb p\equiv
5,7,11  \\ \limfunc{mod}12  \endSb \left( 1+\dsum\limits_{r=1}^\infty \frac
1{p^{2rs}}\right) \cdot \dprod\limits\Sb p\equiv 1  \\ \limfunc{mod}12 
\endSb \left( 1+\dsum\limits_{r=1}^\infty \frac{(r+1)^2}{p^{rs}}\right) \\
\ &=&\dprod\limits\Sb p\equiv 5,7,11  \\ \limfunc{mod}12  \endSb \left(
1-\frac 1{p^s+1}\right) \left( 1-\frac 1{p^s}\right) ^{-1}\cdot
\dprod\limits \Sb p\equiv 1  \\ \limfunc{mod}12  \endSb \left( 1+\frac
2{p^s-1}\right) \left( 1-\frac 1{p^s}\right) ^{-2} \\
\ &=&\left( 1-\frac 1{2^s}\right) \left( 1-\frac 1{3^s}\right) \dprod\limits 
\Sb p\equiv 5,7,11  \\ \limfunc{mod}12  \endSb \left( 1+\frac 1{p^s}\right)
^{-1}\cdot \dprod\limits\Sb p\equiv 1  \\ \limfunc{mod}12  \endSb \left(
1+\frac 1{p^s}\right) \left( 1-\frac 1{p^s}\right) ^{-2}\cdot \zeta (s)
\end{eqnarray*}
and hence 
\[
\dsum\limits_{n\leq N}a_n^2\sim D\,N 
\]
where \cite{Sebah} 
\begin{eqnarray*}
D &=&4\dprod\limits\Sb p\equiv 5,7,11  \\ \limfunc{mod}12  \endSb \left(
1+\frac 1p\right) ^{-1}\cdot \dprod\limits\Sb p\equiv 1  \\ \limfunc{mod}12 
\endSb \left( 1+\frac 1p\right) \left( 1-\frac 1p\right) ^{-2} \\
\ &=&2.6339157938496334172500926....
\end{eqnarray*}
Also, 
\[
\ln \left( \max\limits_{n\leq N}\,|a_n|\right) \sim \ln (2)\frac{\ln (N)}{%
\ln (\ln (N))} 
\]
but a more precise asymptotic statement is evidently open \cite{Tnbm1}.
Finally, 
\[
\dsum\limits_{n\leq N}a_n\sim O\left( N\exp \left( -c(\ln N)^{3/5}(\ln \ln
N)^{-1/5}\right) \right) 
\]
and a stronger conjectured bound $\tsum\nolimits_{n\leq N}a_n=O(N^d)$ with $%
d<1$ is equivalent to a weak form of the generalized Riemann hypothesis \cite
{Tnbm2}.

\subsection{Variation One}

Let the coefficients $\{b_n\}_{n=0}^\infty $ satisfy \cite{Sms1, Sec21} 
\[
\dsum\limits_{n=0}^\infty b_nq^n=\dprod\limits_{m=1}^\infty
(1-q^m)(1-q^{2m}). 
\]
It can be shown that $b_n=g(8n+1)$, where $g$ is multiplicative and 
\[
g(p^r)=\left\{ 
\begin{array}{lll}
1 &  & \text{if }p=2\text{ and }r=1, \\ 
1 &  & \text{if }p\equiv 5,7\limfunc{mod}8\text{ and }r\text{ is even,} \\ 
(-1)^{r/2} &  & \text{if }p\equiv 3\limfunc{mod}8\text{ and }r\text{ is even,%
} \\ 
r+1 &  & \text{if }p\equiv 1\limfunc{mod}8\text{ and }(-4)^{(p-1)/8}\equiv 1%
\limfunc{mod}p\text{,} \\ 
(-1)^r(r+1) &  & \text{if }p\equiv 1\limfunc{mod}8\text{ and }%
(-4)^{(p-1)/8}\equiv -1\limfunc{mod}p\text{,} \\ 
0 &  & \text{otherwise.}
\end{array}
\right. 
\]
For reasons of space, we mention only the mean-square result: 
\[
\dsum\limits_{n\leq N}b_n^2\sim E\,N 
\]
where \cite{Sebah} 
\begin{eqnarray*}
E &=&4\dprod\limits\Sb p\equiv 3,5,7  \\ \limfunc{mod}8  \endSb \left(
1+\frac 1p\right) ^{-1}\cdot \dprod\limits\Sb p\equiv 1  \\ \limfunc{mod}8 
\endSb \left( 1+\frac 1p\right) \left( 1-\frac 1p\right) ^{-2} \\
\ &=&1.7627471740390860504652186....
\end{eqnarray*}

\subsection{Variation Two}

Let the coefficients $\{c_n\}_{n=0}^\infty $ satisfy \cite{Sms1, Sec22} 
\[
\dsum\limits_{n=0}^\infty c_nq^n=\dprod\limits_{m=1}^\infty
(1-q^m)(1-q^{3m}). 
\]
It can be shown that $c_n=h(6n+1)$, where $h$ is multiplicative and 
\[
h(p^r)=\left\{ 
\begin{array}{lll}
1 &  & \text{if }p\equiv 5\limfunc{mod}6\text{ and }r\text{ is even,} \\ 
r+1 &  & \text{if }p\equiv 1\limfunc{mod}6\text{ and }2^{(p-1)/3}\equiv 1%
\limfunc{mod}p\text{,} \\ 
1 &  & \text{if }p\equiv 1\limfunc{mod}6,\text{ }2^{(p-1)/3}\not \equiv 1%
\limfunc{mod}p\text{ and }r\equiv 0\limfunc{mod}3\text{,} \\ 
-1 &  & \text{if }p\equiv 1\limfunc{mod}6,\text{ }2^{(p-1)/3}\not \equiv 1%
\limfunc{mod}p\text{ and }r\equiv 1\limfunc{mod}3\text{,} \\ 
0 &  & \text{otherwise.}
\end{array}
\right. 
\]
Observe that $|h(p)|\,$ assumes two distinct values over the set $p\equiv 1%
\limfunc{mod}6$, which is more complicated than $|f(p)|$ over $p\equiv 1%
\limfunc{mod}12$ and $|g(p)|$ over $p\equiv 1\limfunc{mod}8$. Consequently,
the asymptotics for $\tsum\nolimits_{n\leq N}c_n^2$ are more subtle than
those for $\tsum\nolimits_{n\leq N}a_n^2$ and $\tsum\nolimits_{n\leq N}b_n^2$%
.

\section{Third Power}

Let the coefficients $\{a_n\}_{n=0}^\infty $ satisfy \cite{Ewell, Sec3} 
\[
\dsum\limits_{n=0}^\infty a_nq^n=\dprod\limits_{m=1}^\infty (1-q^m)^3. 
\]
Jacobi's triple-product identity gives that 
\[
a_n=\left\{ 
\begin{array}{lll}
(-1)^k(2k+1) &  & \text{if }8n+1=(2k+1)^2\text{ for some }0\leq k<\infty ,
\\ 
0 &  & \text{otherwise}
\end{array}
\right. 
\]
and hence, by elementary considerations, 
\[
\begin{array}{ccccc}
\dsum\limits_{n\leq N}|a_n|\sim 2N, &  & \dsum\limits_{n\leq N}a_n^2\sim
\dfrac 16(8N)^{3/2}, &  & \max\limits_{n\leq N}\,|a_n|\sim (8N)^{1/2}
\end{array}
\]
as $N\rightarrow \infty $. Also $a_n=f(8n+1)$, where $f$ is multiplicative
and 
\[
f(p^r)=\left\{ 
\begin{array}{lll}
p^{r/2} &  & \text{if }p\equiv 1\limfunc{mod}4\text{ and }r\text{ is even,}
\\ 
(-1)^{r/2}p^{r/2} &  & \text{if }p\equiv 3\limfunc{mod}4\text{ and }r\text{
is even,} \\ 
0 &  & \text{otherwise.}
\end{array}
\right. 
\]
Of course, 
\[
\begin{array}{ccc}
\dsum\limits_{n\leq N}a_n=(-1)^k(k+1) &  & \text{if }(2k+1)^2\leq
8N+1<(2k+3)^2
\end{array}
\]
and thus $\tsum\nolimits_{n\leq N}a_n$ diverges because 
\[
\begin{array}{ccc}
\limfunc{limsup}\limits_{N\rightarrow \infty }\dfrac 1{\sqrt{2N}%
}\dsum\limits_{n\leq N}a_n=1, &  & \limfunc{liminf}\limits_{N\rightarrow
\infty }\dfrac 1{\sqrt{2N}}\dsum\limits_{n\leq N}a_n=-1.
\end{array}
\]

\section{Fourth Power}

Let the coefficients $\{a_n\}_{n=0}^\infty $ satisfy \cite{Rama, Sec4} 
\[
\dsum\limits_{n=0}^\infty a_nq^n=\dprod\limits_{m=1}^\infty (1-q^m)^4. 
\]
No simple formula for $a_n$ is known, but $a_n=f(6n+1)$, where $f$ is
multiplicative and 
\[
f(p^r)=\left\{ 
\begin{array}{lll}
(-1)^{r/2}p^{r/2} &  & \text{if }p\equiv 5\limfunc{mod}6\text{ and }r\text{
is even,} \\ 
\delta _{p,r}\dfrac{(x_p+\sqrt{3}i\,y_p)^{r+1}-(x_p-\sqrt{3}i\,y_p)^{r+1}}{2%
\sqrt{3}i\,y_p} &  & \text{if }p\equiv 1\limfunc{mod}6\text{,} \\ 
0 &  & \text{otherwise}
\end{array}
\right. 
\]
where $i$ is the imaginary unit and $(x_p,y_p)$ is the unique pair of
positive integers for which $p=x_p^2+3y_p^2$. Also, $\delta _{p,r}=-1$ when $%
r$ is odd and $x_p\not \equiv 1\limfunc{mod}3$; otherwise $\delta _{p,r}=1.$

It turns out that $\sum_{n=1}^\infty f(n)n^{-s}$ is the L-series for the
elliptic curve $36A1$: 
\[
y^2=x^3+1 
\]
and hence a theorem of Rankin \cite{Rnkn} implies that 
\[
\dsum\limits_{n\leq N}a_n^2\sim C\,N^2. 
\]
We examine 
\begin{eqnarray*}
\dsum\limits_{n=1}^\infty \frac{f(n)^2}{n^{s+1}} &=&\dprod\limits\Sb p\equiv
5  \\ \limfunc{mod}6  \endSb \left( 1+\dsum\limits_{r=1}^\infty \frac{p^r}{%
p^{2r(s+1)}}\right) \\
&&\ \ \ \ \ \cdot \dprod\limits\Sb p\equiv 1  \\ \limfunc{mod}6  \endSb %
\left( 1-\frac 1{12y_p^2}\dsum\limits_{r=1}^\infty \frac{\left[ \left( x_p+%
\sqrt{3}i\,y_p\right) ^{r+1}-\left( x_p-\sqrt{3}i\,y_p\right) ^{r+1}\right]
^2}{p^{r(s+1)}}\right) \\
\ &=&\dprod\limits\Sb p\equiv 5  \\ \limfunc{mod}6  \endSb \left( 1-\frac
1{p^{2s+1}}\right) ^{-1} \\
&&\ \ \ \ \ \cdot \dprod\limits\Sb p\equiv 1  \\ \limfunc{mod}6  \endSb 
\frac{p^{2(s+1)}\left( p^{s+1}+x_p^2+3y_p^2\right) }{\left[ p^{s+1}-\left(
x_p^2+3y_p^2\right) \right] \left[ \left( x_p^2+3y_p^2\right)
^2-2p^{s+1}\left( x_p^2-3y_p^2\right) +p^{2(s+1)}\right] } \\
\ &=&\dprod\limits\Sb p\equiv 5  \\ \limfunc{mod}6  \endSb \left( 1-\frac
1{p^{2s+1}}\right) ^{-1} \\
&&\ \ \ \ \ \cdot \dprod\limits\Sb p\equiv 1  \\ \limfunc{mod}6  \endSb 
\frac{p^{2(s+1)}\left( p^{s+1}+p\right) }{\left[ p^{s+1}-p\right] \left[
p^2-2p^{s+1}\left( x_p^2-3y_p^2\right) +p^{2(s+1)}\right] } \\
\ &=&\left( 1-\frac 1{2^s}\right) \left( 1-\frac 1{3^s}\right) \dprod\limits 
\Sb p\equiv 5  \\ \limfunc{mod}6  \endSb \left( 1-\frac 1{p^s}\right) \left(
1-\frac 1{p^{2s+1}}\right) ^{-1} \\
&&\ \ \ \ \ \cdot \dprod\limits\Sb p\equiv 1  \\ \limfunc{mod}6  \endSb %
\left( 1+\frac 1{p^s}\right) \left( 1+\frac{2p^{s-1}\left(
x_p^2-3y_p^2\right) -1}{p^{2s}-2p^{s-1}\left( x_p^2-3y_p^2\right) +1}\right)
\cdot \zeta (s)
\end{eqnarray*}
and thus it \emph{seems} that 
\[
C=6\dprod\limits\Sb p\equiv 5  \\ \limfunc{mod}6  \endSb \left( 1-\frac
1p\right) \left( 1-\frac 1{p^3}\right) ^{-1}\cdot \dprod\limits\Sb p\equiv 1 
\\ \limfunc{mod}6  \endSb \left( 1+\frac 1p\right) \left( 1+\frac{2\left(
x_p^2-3y_p^2\right) -1}{p^2-2\left( x_p^2-3y_p^2\right) +1}\right) . 
\]
The language of modular forms is now unavoidable: let 
\[
\begin{array}{ccc}
\eta (t)=e^{\pi it/12}\dprod\limits_{k=1}^\infty \left( 1-e^{2\pi
kit}\right) , &  & \limfunc{Im}(t)>0
\end{array}
\]
denote the Dedekind eta function. The unique cusp form of weight $3$, level $%
12$ and Nebentypus character $(-3/\cdot )$ is \cite{Schuett} 
\[
\begin{array}{ccc}
\eta (2t)^3\eta (6t)^3=q\dprod\limits_{m=1}^\infty \left( 1-q^{2m}\right)
^3\dprod\limits_{m=1}^\infty \left( 1-q^{6m}\right) ^3, &  & q=e^{2\pi it}
\end{array}
\]
and possesses the expansion 
\[
L_{xy}(s)=\left( 1+\frac 1{3^{s-1}}\right) ^{-1}\dprod\limits\Sb p\equiv 5 
\\ \limfunc{mod}6  \endSb \left( 1-\frac 1{p^{2(s-1)}}\right) ^{-1}\cdot
\dprod\limits\Sb p\equiv 1  \\ \limfunc{mod}6  \endSb \left( 1-\frac{2\left(
x_p^2-3y_p^2\right) }{p^s}+\frac 1{p^{2(s-1)}}\right) ^{-1}; 
\]
hence 
\begin{eqnarray*}
L_{xy}(2) &=&\frac 34\dprod\limits\Sb p\equiv 5  \\ \limfunc{mod}6  \endSb %
\left( 1-\frac 1{p^2}\right) ^{-1}\cdot \dprod\limits\Sb p\equiv 1  \\ 
\limfunc{mod}6  \endSb \left( 1+\frac{2\left( x_p^2-3y_p^2\right) -1}{%
p^2-2\left( x_p^2-3y_p^2\right) +1}\right) \\
\ &=&0.7372929961855962401764261...;
\end{eqnarray*}
hence \cite{Sebah} 
\begin{eqnarray*}
C &=&8\,L_{xy}(2)\dprod\limits\Sb p\equiv 5  \\ \limfunc{mod}6  \endSb %
\left( 1-\frac 1p\right) \left( 1-\frac 1{p^2}\right) \left( 1-\frac
1{p^3}\right) ^{-1}\cdot \dprod\limits\Sb p\equiv 1  \\ \limfunc{mod}6 
\endSb \left( 1+\frac 1p\right) \\
\ &=&4.6417183001293981350615666....
\end{eqnarray*}
Observe, however, that in the vicinity of $s=1$, the series 
\[
\dsum\limits_p\frac{2p^{s-1}\left( x_p^2-3y_p^2\right) -1}{%
p^{2s}-2p^{s-1}\left( x_p^2-3y_p^2\right) +1} 
\]
is not absolutely convergent since either $x_p^2$ or $3y_p^2$ is often of
size $p$. The requirement that $(\sum f(n)^2n^{-s-1})\zeta (s)^{-1}$ is
analytic consequently might not hold. Empirically, we compute 
\[
\frac 1{N^2}\dsum\limits_{n\leq N}a_n^2\approx 4.877... 
\]
which raises further doubt about the correctness of our value $4.641...$ for 
$C.$

\subsection{Variation One}

Let the coefficients $\{b_n\}_{n=0}^\infty $ satisfy \cite{Gls2, Sec41} 
\[
\dsum\limits_{n=0}^\infty b_nq^n=\dprod\limits_{m=1}^\infty
(1-q^m)^2(1-q^{2m})^2. 
\]
It can be shown that $b_n=g(4n+1)$, where $g$ is multiplicative and 
\[
g(p^r)=\left\{ 
\begin{array}{lll}
(-1)^{r/2}p^{r/2} &  & \text{if }p\equiv 3\limfunc{mod}4\text{ and }r\text{
is even,} \\ 
\varepsilon _{p,r}\dfrac{(u_p+i\,v_p)^{r+1}-(u_p-i\,v_p)^{r+1}}{2i\,v_p} & 
& \text{if }p\equiv 1\limfunc{mod}4\text{,} \\ 
0 &  & \text{otherwise}
\end{array}
\right. 
\]
where $(u_p,v_p)$ is the unique pair of positive integers for which $%
p=u_p^2+v_p^2$, $u_p$ is odd and $v_p$ is even. Also, $\varepsilon
_{p,r}=(-1)^{(u_p+v_p-1)/2}$ when $r$ is odd and $\varepsilon _{p,r}=1$ when 
$r$ is even.

It turns out that $\sum_{n=1}^\infty g(n)n^{-s}$ is the L-series for the
elliptic curve $32A2$: 
\[
y^2=x^3-x 
\]
and hence Rankin's theorem implies that 
\[
\dsum\limits_{n\leq N}b_n^2\sim D\,N^2 
\]
where \emph{seemingly} 
\[
D=4\dprod\limits\Sb p\equiv 3  \\ \limfunc{mod}4  \endSb \left( 1-\frac
1p\right) \left( 1-\frac 1{p^3}\right) ^{-1}\cdot \dprod\limits\Sb p\equiv 1 
\\ \limfunc{mod}4  \endSb \left( 1+\frac 1p\right) \left( 1+\frac{2\left(
u_p^2-v_p^2\right) -1}{p^2-2\left( u_p^2-v_p^2\right) +1}\right) . 
\]
The unique cusp form of weight $3$, level $16$ and Nebentypus character $%
(-4/\cdot )$ is \cite{Schuett} 
\[
\eta (4t)^6=q\dprod\limits_{m=1}^\infty \left( 1-q^{4m}\right) ^6 
\]
and possesses the expansion 
\[
L_{uv}(s)=\dprod\limits\Sb p\equiv 3  \\ \limfunc{mod}4  \endSb \left(
1-\frac 1{p^{2(s-1)}}\right) ^{-1}\cdot \dprod\limits\Sb p\equiv 1  \\ 
\limfunc{mod}4  \endSb \left( 1-\frac{2\left( u_p^2-v_p^2\right) }{p^s}%
+\frac 1{p^{2(s-1)}}\right) ^{-1}; 
\]
hence 
\begin{eqnarray*}
L_{uv}(2) &=&\dprod\limits\Sb p\equiv 3  \\ \limfunc{mod}4  \endSb \left(
1-\frac 1{p^2}\right) ^{-1}\cdot \dprod\limits\Sb p\equiv 1  \\ \limfunc{mod}%
4  \endSb \left( 1+\frac{2\left( u_p^2-v_p^2\right) -1}{p^2-2\left(
u_p^2-v_p^2\right) +1}\right) \\
\ &=&0.8593982272525466034362619...;
\end{eqnarray*}
hence \cite{Sebah} 
\begin{eqnarray*}
D &=&4\,L_{uv}(2)\dprod\limits\Sb p\equiv 3  \\ \limfunc{mod}4  \endSb %
\left( 1-\frac 1p\right) \left( 1-\frac 1{p^2}\right) \left( 1-\frac
1{p^3}\right) ^{-1}\cdot \dprod\limits\Sb p\equiv 1  \\ \limfunc{mod}4 
\endSb \left( 1+\frac 1p\right) \\
\ &=&1.9533514553987911733090376....
\end{eqnarray*}
Again, analyticity requirements might not hold and we empirically compute 
\[
\frac 1{N^2}\dsum\limits_{n\leq N}b_n^2\approx 2.188... 
\]
which raises doubt about our value $1.953...$ for $D.$

\subsection{Variation Two}

Let the coefficients $\{c_n\}_{n=0}^\infty $ satisfy \cite{Sms2, Sec42} 
\[
\dsum\limits_{n=0}^\infty c_nq^n=\dprod\limits_{m=1}^\infty
(1-q^m)^2(1-q^{3m})^2. 
\]
It can be shown that $c_n=h(3n+1)$, where $h$ is multiplicative and 
\[
h(p^r)=\left\{ 
\begin{array}{lll}
(-1)^{r/2}p^{r/2} &  & \text{if }p\equiv 2\limfunc{mod}3\text{ and }r\text{
is even,} \\ 
\tilde \delta _{p,r}\dfrac{\left( \dfrac{z_p+3\sqrt{3}i\,w_p}2\right)
^{r+1}-\left( \dfrac{z_p-3\sqrt{3}i\,w_p}2\right) ^{r+1}}{3\sqrt{3}i\,w_p} & 
& \text{if }p\equiv 1\limfunc{mod}3\text{,} \\ 
0 &  & \text{otherwise}
\end{array}
\right. 
\]
where $(z_p,w_p)$ is the unique pair of positive integers for which $%
4p=z_p^2+27w_p^2$. Also, $\tilde \delta _{p,r}=-1$ when $r$ is odd and $%
x_p\equiv 1\limfunc{mod}3$; otherwise $\tilde \delta _{p,r}=1.$

It turns out that $\sum_{n=1}^\infty h(n)n^{-s}$ is the L-series for the
elliptic curve $27A3$: 
\[
y^2+y=x^3 
\]
and hence Rankin's theorem implies that 
\[
\dsum\limits_{n\leq N}c_n^2\sim E\,N^2 
\]
where \emph{seemingly} 
\[
E=3\dprod\limits\Sb p\equiv 2  \\ \limfunc{mod}3  \endSb \left( 1-\frac
1p\right) \left( 1-\frac 1{p^3}\right) ^{-1}\cdot \dprod\limits\Sb p\equiv 1 
\\ \limfunc{mod}3  \endSb \left( 1+\frac 1p\right) \left( 1+\frac{%
(z_p^2-27w_p^2)/2-1}{p^2-(z_p^2-27w_p^2)/2+1}\right) . 
\]
The vector space of cusp forms of weight $3$, level $27$ and Nebentypus
character $(-3/\cdot )$ is three-dimensional. A certain basis element is
given by \cite{Schuett, Sms0} 
\[
\eta (9t)\,\eta (3t)^2\left( \eta (3t)^3+9\,\eta (27t)^3\right) 
\]
and possesses the expansion 
\[
L_{zw}(s)=\dprod\limits\Sb p\equiv 2  \\ \limfunc{mod}3  \endSb \left(
1-\frac 1{p^{2(s-1)}}\right) ^{-1}\cdot \dprod\limits\Sb p\equiv 1  \\ 
\limfunc{mod}3  \endSb \left( 1-\frac{z_p^2-27w_p^2}{2\,p^s}+\frac
1{p^{2(s-1)}}\right) ^{-1}; 
\]
hence 
\begin{eqnarray*}
L_{zw}(2) &=&\dprod\limits\Sb p\equiv 2  \\ \limfunc{mod}3  \endSb \left(
1-\frac 1{p^2}\right) ^{-1}\cdot \dprod\limits\Sb p\equiv 1  \\ \limfunc{mod}%
3  \endSb \left( 1+\frac{(z_p^2-27w_p^2)/2-1}{p^2-(z_p^2-27w_p^2)/2+1}\right)
\\
\ &=&1.0403374913367121372113004...;
\end{eqnarray*}
hence \cite{Sebah} 
\begin{eqnarray*}
E &=&3\,L_{zw}(2)\dprod\limits\Sb p\equiv 2  \\ \limfunc{mod}3  \endSb %
\left( 1-\frac 1p\right) \left( 1-\frac 1{p^2}\right) \left( 1-\frac
1{p^3}\right) ^{-1}\cdot \dprod\limits\Sb p\equiv 1  \\ \limfunc{mod}3 
\endSb \left( 1+\frac 1p\right) \\
\ &=&1.0526097875093498936749762....
\end{eqnarray*}
Again, analyticity requirements might not hold and we empirically compute 
\[
\frac 1{N^2}\dsum\limits_{n\leq N}c_n^2\approx 1.290... 
\]
which raises doubt about our value $1.052...$ for $E.$

\section{Sixth Power}

Let the coefficients $\{a_n\}_{n=0}^\infty $ satisfy \cite{Rama, Sec5} 
\[
\dsum\limits_{n=0}^\infty a_nq^n=\dprod\limits_{m=1}^\infty (1-q^m)^6. 
\]
No simple formula for $a_n$ is known, but $a_n=f(4n+1)$, where $f$ is
multiplicative and 
\[
f(p^r)=\left\{ 
\begin{array}{lll}
p^r &  & 
\begin{array}{c}
\text{if }p\equiv 3\limfunc{mod}4\text{ } \\ 
\text{and }r\text{ is even,}
\end{array}
\\ 
(-1)^r\dfrac{(v_p^2-u_p^2+2i\,u_p\,v_p)^{r+1}-(v_p^2-u_p^2-2i\,u_p%
\,v_p)^{r+1}}{4i\,u_pv_p} &  & 
\begin{array}{c}
\text{if }p\equiv 1\limfunc{mod}4\text{,}
\end{array}
\\ 
0 &  & 
\begin{array}{c}
\text{otherwise}
\end{array}
\end{array}
\right. 
\]
where $(u_p,v_p)$ is the unique pair of positive integers for which $%
p=u_p^2+v_p^2$, $u_p$ is odd and $v_p$ is even (as before). Rankin's theorem
implies that 
\[
\dsum\limits_{n\leq N}a_n^2\sim C\,N^3 
\]
and the Selberg-Delange method perhaps can be used to compute $C$.

\section{Eighth Power}

Let the coefficients $\{a_n\}_{n=0}^\infty $ satisfy \cite{Rama, Sec6} 
\[
\dsum\limits_{n=0}^\infty a_nq^n=\dprod\limits_{m=1}^\infty (1-q^m)^8. 
\]
No simple formula for $a_n$ is known, but $a_n=f(3n+1)$, where $f$ is
multiplicative and 
\[
f(p^r)=\left\{ 
\begin{array}{lll}
(-1)^{r/2}p^{3r/2} &  & 
\begin{array}{c}
\text{if }p\equiv 2\limfunc{mod}3\text{ } \\ 
\text{and }r\text{ is even,}
\end{array}
\\ 
\delta _{p,r}\tfrac{\left[ x_p(x_p^2-9y_p^2)+3\sqrt{3}i\,y_p(x_p^2-y_p^2)%
\right] ^{r+1}-\left[ x_p(x_p^2-9y_p^2)-3\sqrt{3}i\,y_p(x_p^2-y_p^2)\right]
^{r+1}}{6\sqrt{3}i\,y_p(x_p^2-y_p^2)} &  & 
\begin{array}{c}
\text{if }p\equiv 1\limfunc{mod}6\text{,}
\end{array}
\\ 
0 &  & 
\begin{array}{c}
\text{otherwise}
\end{array}
\end{array}
\right. 
\]
where $(x_p,y_p)$ is the unique pair of positive integers for which $%
p=x_p^2+3y_p^2$. Also, $\delta _{p,r}=-1$ when $r$ is odd and $x_p\not
\equiv 1\limfunc{mod}3$; otherwise $\delta _{p,r}=1$ (as before). Rankin's
theorem implies that 
\[
\dsum\limits_{n\leq N}a_n^2\sim C\,N^4 
\]
and again the Selberg-Delange method perhaps can be used to compute $C$.

\section{Tenth Power}

For divisors $\ell $ of $24$ (as in the preceding sections), the
multiplicative function $f$ satisfying 
\[
\begin{array}{ccc}
\dsum\limits_{n=0}^\infty a_nq^n=\dprod\limits_{m=1}^\infty (1-q^m)^\ell , & 
& a_n=f(24n/\ell +1)
\end{array}
\]
obeys the standard recurrence 
\[
\begin{array}{ccc}
f(p^r)-f(p)f(p^{r-1})+p^{(\ell -2)/2}f(p^{r-2})=0, &  & r\geq 2.
\end{array}
\]
This is also true when $\ell =10$ except we have $a_n=f(12n+5)/48$.

For the case $p\equiv 5\limfunc{mod}12$, the required initial condition is 
\cite{Sec7, Wnqst, BCLY, Liu, Chu} 
\[
f(p)=8\varphi (u_p,v_p)u_pv_p\left( u_p^2-v_p^2\right) 
\]
where $(u_p,v_p)$ is as before and 
\[
\varphi (u,v)=\left\{ 
\begin{array}{lll}
-1 &  & \text{if }\left( u\equiv 1\limfunc{mod}6\text{ and }v\not \equiv 4%
\limfunc{mod}6\right) \text{ or }\left( u\not \equiv 1\limfunc{mod}6\text{
and }v\equiv 4\limfunc{mod}6\right) , \\ 
1 &  & \text{otherwise.}
\end{array}
\right. 
\]
For the case $p\equiv 1\limfunc{mod}12$, the required initial condition is 
\cite{Sec7, Sms3}: 
\[
f(p)=2\theta (u_p)\left( 2u_pv_p+u_p^2-v_p^2\right) \left(
2u_pv_p+v_p^2-u_p^2\right) 
\]
where $\theta (u)=1$ when $u\equiv 3\limfunc{mod}6$; otherwise $\theta
(u)=-1 $. We summarize the remaining cases: 
\[
f(p^r)=\left\{ 
\begin{array}{lll}
p^{2r} &  & \text{if }p\equiv 7,11\limfunc{mod}12\text{ and }r\text{ is even,%
} \\ 
0 &  & \text{if }\left( p\equiv 7,11\limfunc{mod}12\text{ and }r\text{ is odd%
}\right) \text{ or }\left( p=2,3\right) \text{.}
\end{array}
\right. 
\]

\section{Fourteenth Power}

The standard recurrence is only partly true when $\ell =14$: 
\[
\begin{array}{ccc}
f(p^r)-f(p)f(p^{r-1})\pm p^{(\ell -2)/2}f(p^{r-2})=0, &  & r\geq 2
\end{array}
\]
where the $-$ symbol is chosen when $p\equiv 7\limfunc{mod}12$ and the $+$
symbol is chosen otherwise. We have $a_n=f(12n+7)/(360\sqrt{-3})$; it is
interesting that complex algebraic integers enter the formulation here.

For the case $p\equiv 7\limfunc{mod}12$, the required initial condition is 
\cite{Wnqst, Sec8, CHL} 
\[
f(p)=12(-1)^{(p-7)/12}\psi (x_p,y_p)x_py_p\left( x_p^2-y_p^2\right) \left(
9y_p^2-x_p^2\right) \sqrt{-3} 
\]
where $(x_p,y_p)$ is as earlier and 
\[
\psi (x,y)=\left\{ 
\begin{array}{lll}
-1 &  & \text{if }\left( x\equiv 2\limfunc{mod}6\text{ and }y\not \equiv 1%
\limfunc{mod}4\right) \text{ or }\left( x\not \equiv 2\limfunc{mod}6\text{
and }y\equiv 1\limfunc{mod}4\right) , \\ 
1 &  & \text{otherwise.}
\end{array}
\right. 
\]
For the case $p\equiv 1\limfunc{mod}12$, the required initial condition is 
\cite{Sec8, Sms4}: 
\[
f(p)=2\omega (x_p,y_p)\left( x_p^2-3y_p^2\right) \left(
6x_py_p+x_p^2-3y_p^2\right) \left( 6x_py_p+3y_p^2-x_p^2\right) 
\]
where $\omega (x,y)=(-1)^{y/2}$ when $x\equiv 5,7\limfunc{mod}12$; otherwise 
$\omega (x,y)=(-1)^{y/2+1}$. We summarize the remaining cases: 
\[
f(p^r)=\left\{ 
\begin{array}{lll}
(-1)^{r/2}p^{3r} &  & \text{if }p\equiv 5\limfunc{mod}12\text{ and }r\text{
is even,} \\ 
p^{3r} &  & \text{if }p\equiv 11\limfunc{mod}12\text{ and }r\text{ is even,}
\\ 
0 &  & \text{if }\left( p\equiv 5,11\limfunc{mod}12\text{ and }r\text{ is odd%
}\right) \text{ or }\left( p=2,3\right) \text{.}
\end{array}
\right. 
\]

\section{Twenty-Sixth Power}

For $\ell =26$, we merely point out that \cite{Sec9, Dysn, CCT1, CCT2} 
\begin{eqnarray*}
a_{(p-13)/12} &=&\frac{(-1)^{\kappa _p}}{1019304}\left[ \left(
u_p^2-x_p^2\right) \left( v_p^2-x_p^2\right) \left( \left( 2u_p+x_p\right)
^2-9y_p^2\right) \cdot \right. \\
&&\ \ \left. \left( \left( 2u_p-x_p\right) ^2-9y_p^2\right) \left( \left(
2v_p+x_p\right) ^2-9y_p^2\right) \left( \left( 2v_p-x_p\right)
^2-9y_p^2\right) \right]
\end{eqnarray*}
for primes $p$ of the form $p\equiv 1\limfunc{mod}12$, where $(u_p,v_p)$, $%
(x_p,y_p)$ are as earlier. The quantity $\kappa _p$ is $u_p$ if $u_p\equiv 0%
\limfunc{mod}3$; otherwise $\kappa _p$ is $v_p$. This is a first step toward
writing $a_n$ in terms of a multiplicative function $f$.

\section{Closing Words}

For the coefficients $\{s_n\}_{n=0}^\infty $ and $\{t_n\}_{n=0}^\infty $
with 
\[
\begin{array}{ccc}
\dsum\limits_{n=0}^\infty s_nq^n=\dprod\limits_{m=1}^\infty (1-q^m)^{12}, & 
& \dsum\limits_{n=0}^\infty t_nq^n=\dprod\limits_{m=1}^\infty (1-q^m)^{24},
\end{array}
\]
we have 
\[
\begin{array}{ccc}
\dsum\limits_{n\leq N}s_n^2\sim A\,N^6, &  & \dsum\limits_{n\leq N}t_n^2\sim
B\,N^{12}
\end{array}
\]
as $N\rightarrow \infty $. Clearly $t_n=\tau (n+1)$, where $\tau $ is the
famous Ramanujan tau function. Let also $s_n=\sigma (2n+1)$ for convenience.
We mention Lehmer's conjectures: 
\[
\begin{array}{ccc}
s_n\neq 0\text{ for all }n, &  & t_n\neq 0\text{ for all }n.
\end{array}
\]
By multiplicativity of $\sigma $ and $\tau $, it would be sufficient to
prove that 
\[
\begin{array}{ccc}
\sigma (p^r)\neq 0\text{ for all }p>2, &  & \tau (p^r)\neq 0\text{ for all }%
p\geq 2
\end{array}
\]
but no formulas for $\sigma (p^r)$ or $\tau (p^r)$ are available (unlike
earlier). It can be demonstrated, however, that $B$ is a product of certain
L-series values and hence \cite{Lhm, Zgr, KZ, F1} 
\[
B=0.0320070045390141974938639.... 
\]
A\ corresponding result for $A$ is open, as far as is known. Information on
numerical estimates of $A$ would be greatly appreciated.

\section{Acknowledgements}

I thank the following generous individuals for their help:

\begin{itemize}
\item  Michael Somos for developing $f(p^r)$ expressions in sections 2 \&\ 4
and for insights leading to sections 7 \&\ 8,

\item  Pascal Sebah for computing high-precision numerical estimates of
infinite products,

\item  Tim Dokchitser for writing ComputeL (which was used to evaluate
L-series at $s=2$),

\item  Matthias Schuett for pointing out the three ``complementary'' cusp
forms in section 4,

\item  G\'erald Tenenbaum for his expertise in the Selberg-Delange method,

\item  Peter Stevenhagen for simplifying $f(p^r)$ expressions in section 2,

\item  Kevin Buzzard for discussing the absolute convergence issue in
section 4.
\end{itemize}

\noindent If anyone makes progress on the unsolved problems mentioned in
this paper, I would be most grateful for a notification.

\end{document}